\RequirePackage{ifpdf}
\ifpdf 
\documentclass[pdftex]{sigma}
\else
\documentclass{sigma}
\fi

\numberwithin{equation}{section}
\numberwithin{theorem}{section}
\numberwithin{proposition}{section}
\numberwithin{lemma}{section}
\numberwithin{corollary}{section}
\numberwithin{definition}{section}
\numberwithin{example}{section}
\numberwithin{remark}{section}
\numberwithin{note}{section}

\newcommand{\R}{\mathbb R}
\newcommand{\N}{\mathbb N}
\newcommand{\C}{\mathbb C}
\newcommand{\Z}{\mathbb Z}
\newcommand{\T}{\mathbb T}
\newcommand{\eps}{\varepsilon}
\newcommand{\dirint}[2]{\sideset{}{^\oplus}\int\limits_{#1\ }^{\ \ #2}}

\newcommand{\rphis}[5]{\,_{#1}\varphi_{#2} \left( \genfrac{.}{.}{0pt}{}{#3}{#4}
\ ;#5 \right)}
\newcommand{\su}{\mathfrak{su}}
\newcommand{\U}{\mathcal U}
\newcommand{\Res}[1]{\underset{#1}{\mathrm{Res}}}
\newcommand{\vect}[2]{\left( \genfrac{.}{.}{0pt}{}{#1}{#2} \right) }

\begin{document}

\allowdisplaybreaks

\renewcommand{\thefootnote}{$\star$}

\renewcommand{\PaperNumber}{077}

\FirstPageHeading

\ShortArticleName{Quantum Analogs of Tensor Product Representations of $\su(1,1)$}

\ArticleName{Quantum Analogs\\ of Tensor Product Representations of $\boldsymbol{\su(1,1)}\,$\footnote{This paper is a
contribution to the Special Issue ``Relationship of Orthogonal Polynomials and Special Functions with Quantum Groups and Integrable Systems''. The
full collection is available at
\href{http://www.emis.de/journals/SIGMA/OPSF.html}{http://www.emis.de/journals/SIGMA/OPSF.html}}}

\Author{Wolter GROENEVELT}

\AuthorNameForHeading{W.~Groenevelt}

\Address{Delft Institute of Applied Mathematics, Technische Universiteit Delft,\\ PO Box 5031,
2600 GA Delft, the Netherlands}
\Email{\href{mailto:w.g.m.groenevelt@tudelft.nl}{w.g.m.groenevelt@tudelft.nl}}
\URLaddress{\url{http://fa.its.tudelft.nl/~groenevelt/}}

\ArticleDates{Received April 28, 2011, in f\/inal form August 04, 2011;  Published online August 09, 2011}

\Abstract{We study representations of $\U_q(\su(1,1))$ that can be considered as quantum analogs of tensor products of irreducible $*$-representations of the Lie algebra $\su(1,1)$. We determine the decomposition of these representations into irreducible $*$-representations of $\U_q(\su(1,1))$ by diagonalizing the action of the Casimir operator on suitable subspaces of the representation spaces. This leads to an interpretation of the big $q$-Jacobi polynomials and big $q$-Jacobi functions as quantum analogs of Clebsch--Gordan coef\/f\/icients.}


\Classification{20G42; 33D80}

\section{Introduction}

The quantum algebra $\U_q=\U_q(\su(1,1))$ has f\/ive classes of irreducible $*$-representations: the po\-si\-ti\-ve and negative discrete series~$\pi^\pm$, the principal unitary series~$\pi^{\mathrm{P}}$, the complementary series~$\pi^{\mathrm{C}}$, and the strange series~$\pi^{\mathrm{S}}$. The f\/irst four classes of representations can be considered ``classical'' in the sense that they are natural quantum analogs of the four irreducible $*$-representations of the Lie algebra $\su(1,1)$. In the classical limit $q \uparrow 1$ these representations all tend to their classical counterparts. The f\/ifth class has no classical analog, hence the name ``strange series''. This class of representations disappears in the classical limit.

In this paper we study representations of $\U_q$ that can be considered as quantum analogs of tensor products of irreducible $*$-representations of $\su(1,1)$, but the representations that we consider are not tensor products of irreducible $*$-representations of $\U_q$. The motivation for studying such representations comes from corepresentation theory of the locally compact quantum group analog $M$ of the normalizer of $SU(1,1)$ in $SL(2,\C)$. The dual quantum group $\widehat M$ is generated as a von Neumann algebra by the standard generators of $\U_q$ and two extra generators. In this sense $\U_q$ can be considered as a subalgebra of $\widehat M$.  An irreducible discrete series representation of $\widehat M$ restricted to $\U_q$ decomposes as the sum $\pi^+ \oplus \pi^- \oplus \pi^{\mathrm{S}}$, with appropriate representation labels, see \cite[Section 5]{GKK}. A tensor product of such representations consists of a sum of nine simple tensor products. For f\/ive of these simple tensor products it is known how to decompose them into irreducible $\U_q$-representations: $\pi^+\otimes \pi^+$, $\pi^- \otimes \pi^-$, $\pi^+\otimes \pi^-$, $\pi^+ \otimes \pi^{\mathrm{S}}$ and $\pi^{\mathrm{S}} \otimes \pi^-$, see e.g.~\cite[Section 4]{KMM}, \cite[Section 2]{KM}, \cite[Section 8]{Gr06}. In this paper we consider the remaining terms as two ``indivisible" representations, $\mathcal T =(\pi^- \otimes \pi^+) \oplus (\pi^{\mathrm{S}} \otimes \pi^{\mathrm{S}})$ and $\mathcal T'=(\pi^- \otimes \pi^{\mathrm{S}}) \oplus (\pi^{\mathrm{S}} \otimes \pi^+)$, and determine their decompositions. In a similar way the principal unitary series and complementary series representations of $\widehat M$ restricted to $\U_q$ decompose as $\pi^{\mathrm{P}} \oplus \pi^{\mathrm{P}}$ and $\pi^{\mathrm{C}} \oplus \pi^{\mathrm{C}}$, respectively. Taking tensor products of these we end up again with ``indivisible" sums of simple tensor products that can be considered as natural analogs of the tensor product of two principal unitary series or complementary series of $\su(1,1)$. It is remarkable that the decomposition of the representations we consider here were already announced in~\cite[Section~12]{KV}.

The Clebsch--Gordan coef\/f\/icients with respect to standard bases for the f\/ive simple tensor products mentioned above can be described in terms of terminating basic hypergeometric $_3\varphi_2$-series. The orthogonality relations of the Clebsch--Gordan coef\/f\/icients correspond to the orthogonality relations of the (dual) $q$-Hahn polynomials and the continuous dual $q$-Hahn polynomials. The Clebsch--Gordan coef\/f\/icients for the representations we consider in this paper turn out to be non-terminating $_3\varphi_2$-series, and consequently the corresponding orthogonality relations are (in general) not related to orthogonal polynomials, but to non-polynomial unitary transform pairs.

Let us now brief\/ly describe the contents of this paper. In Section~\ref{sec:Uq} we recall the def\/inition of the quantum algebra $\U_q(\su(1,1))$ and its irreducible $*$-representations. In Section~\ref{sec:-+} we consider the decomposition of~$\mathcal T$. Our choice of representation labels is slightly more general than allowed in the context of the locally compact quantum group~$M$. The representation $\mathcal T$ can be considered as a quantum analog of the tensor product of a negative and a positive discrete series representation of $\su(1,1)$. We diagonalize the action of the Casimir operator, and this naturally leads to the interpretation of big $q$-Jacobi functions as quantum analogs of Clebsch--Gordan coef\/f\/icients. We also consider the representation $\mathcal T'$, which completes $\mathcal T$ to a genuine tensor product representation of $\U_q$, but has no classical analog. In Section~\ref{sec:PP} we consider the representation $(\pi^{\mathrm{P}} \otimes \pi^{\mathrm{P}}) \oplus (\pi^{\mathrm{P}} \otimes \pi^{\mathrm{P}})$. The diagonalization of the Casimir operator leads in this case to vector-valued big $q$-Jacobi functions as Clebsch--Gordan coef\/f\/icients. Finally in Section~\ref{sec:more} we give the decompositions of several other quantum analogs of tensor product representations.

\textbf{Notations.}
We use $\N=\{0,1,2,\ldots\}$ and $q$ is a f\/ixed number in $(0,1)$. We use standard notations for $q$-shifted factorials, theta functions and basic hypergeometric series from the book of Gasper and Rahman \cite{GRa04}. For $x \in \C$ and $n \in \N \cup \{\infty\}$ the $q$-shifted factorial is def\/ined by $(x;q)_n = \prod\limits_{k=0}^{n-1} (1-xq^{k})$, where the empty product is equal to~$1$. For $x \neq 0$ the normalized Jacobi $\theta$-function is def\/ined by $\theta(x;q)=(x;q)_\infty (q/x;q)_\infty$. For products of $q$-shifted factorial and products of $\theta$-functions we use the notations
\[
(x_1,x_2,\ldots,x_m;q)_n=\prod_{k=1}^m (x_k;q)_n, \qquad  \theta(x_1,x_2,\ldots,x_m;q)= \prod_{k=1}^m \theta(x_k;q)
\]
and
\[
\big(x^{\pm 1};q\big)_n = (x,1/x;q)_n, \qquad   \theta\big(x^{\pm 1};q\big) = \theta(x,1/x;q).
\]
The basic hypergeometric series $_3\varphi_2$ is def\/ined by
\[
\rphis{3}{2}{a,b,c}{d,e}{q,z} = \sum_{k=0}^\infty \frac{ (a,b,c;q)_k }{ (q,d,e;q)_k } z^k.
\]

\section[The quantum algebra $\U_q(\su(1,1))$]{The quantum algebra $\boldsymbol{\U_q(\su(1,1))}$} \label{sec:Uq}

The quantized universal enveloping algebra $\mathcal U_q = \U_q\big(\su(1,1)\big)$ is the unital, associative, complex algebra generated by $K$, $K^{-1}$, $E$, and $F$, subject to the relations
\begin{gather}
K K^{-1} = 1 = K^{-1}K, \qquad KE = qEK, \qquad KF= q^{-1}FK, \nonumber\\
 EF-FE =\frac{K^2-K^{-2}}{q-q^{-1}}.\label{eq:commutationrel}
\end{gather}
The Casimir element
\begin{gather} \label{eq:Casimir}
\Omega = \frac{q^{-1} K^2 +qK^{-2}-2}{(q^{-1}-q)^2} + EF= \frac{q^{-1}K^{-2}+qK^2-2}{(q^{-1}-q)^2} +FE
\end{gather}
is a central element of $\U_q$. The algebra
$\U_q$ is a Hopf $*$-algebra with comultiplication $\Delta$ def\/ined on the generators by
\begin{alignat}{3}
& \Delta(K) = K \otimes K,\qquad && \Delta(E)= K \otimes E + E \otimes K^{-1},& \nonumber\\
& \Delta(K^{-1}) = K^{-1} \otimes K^{-1},\qquad && \Delta(F) = K \otimes F + F \otimes K^{-1} ,& \label{eq:comult}
\end{alignat}
and $\Delta$ is extended to $\mathcal U_q$ as an algebra homomorphism. In particular, it follows from~\eqref{eq:Casimir} and~\eqref{eq:comult} that
\begin{gather}
\Delta(\Omega) = \frac{1}{(q^{-1}-q)^{2}} \big[ q^{-1}\big( K^2 \otimes K^2 \big)+ q
\big(K^{-2} \otimes K^{-2}\big) -2 (1\otimes 1) \big] \nonumber\\
\phantom{\Delta(\Omega) =}{} + K^2 \otimes EF +KF\otimes
EK^{-1} + EK \otimes K^{-1}F + EF \otimes K^{-2}.\label{eq:De(OM)}
\end{gather}
The $*$-structure on $\U_q$ is def\/ined on the generators by
\[
K^*=K, \qquad E^*=-F, \qquad F^* = -E, \qquad (K^{-1})^* = K^{-1}.
\]
Note that the Casimir element is self-adjoint in $\U_q$, i.e.~$\Omega^*=\Omega$.

There are, besides the trivial representation, f\/ive classes of irreducible $*$-representations of~$\U_q$, see \cite[Proposition~4]{MMNNSU}, \cite[Section~6]{BK93}. The representations are given in terms of unbounded operators on $\ell^2(\N)$ or $\ell^2(\Z)$. As common domain we take the f\/inite linear combinations of the standard orthonormal basis vectors $e_n$. The representations are unbounded $*$-representations in the sense of Schm\"udgen \cite[Def\/inition~8.1.9]{Sch}. Below we list the actions of the generators $K$, $K^{-1}$, $E$, $F$ on basis vectors $e_n$. The Casimir element~$\Omega$ plays an important role in this paper, therefore we also list the action of $\Omega$.

The representation listed below are completely characterized, up to unitary equivalence, by the actions of~$K$ and~$\Omega$. Let us brief\/ly describe how these actions determine the actions of~$E$ and~$F$ up to a phase factor. Let $\pi$ be a $\U_q$-representation acting on $\ell^2(\Z)$ with orthonormal basis $\{e_n\}_{n \in \Z}$, and suppose that $\pi(K)e_n=q^{n+\eps}e_n$ and $\pi(\Omega)e_n = \omega e_n$ for all $n \in \Z$, where $\eps$ and $\omega$ are real numbers. The commutation relations~\eqref{eq:commutationrel} imply that $\pi(E)e_n = c_n e_{n+1}$ and $\pi(F) e_n = d_n e_{n-1}$ for certain numbers~$c_n$ and~$d_n$. Furthermore, the relation $E^*=-F$ implies $c_n=-\overline{d_{n+1}}$, so that $\pi(FE) e_n = - |c_n|^2 e_n$. On the other hand, from~\eqref{eq:Casimir} and the actions of $K$ and $\Omega$ it follows that
\[
\pi(FE) e_n = \left(\omega-\frac{q^{2n+2\eps+1}+q^{-2n-2\eps-1}-2}{(q^{-1}-q)^2}\right)e_n.
\]
So $c_n$ can be determined up to a phase factor.

The f\/ive classes of irreducible $*$-representations of $\U_q$ are the following:

\textbf{Positive discrete series.} The positive discrete series $\pi^+_k$
are labeled by $k>0$. The representation space is $\ell^2(\N)$
with orthonormal basis $\{e_n\}_{n \in \N }$. The action is
given by
\begin{gather}
\pi^+_k(K)  e_n  = q^{k+n}  e_n, \qquad \pi^+_k\big(K^{-1}\big)  e_n =q^{-(k+n)} e_n, \nonumber \\
\big(q^{-1}-q\big)\pi^+_k(E)  e_n  = q^{-\frac12-k-n} \sqrt{(1-q^{2n+2})(1-q^{4k+2n})}  e_{n+1},  \nonumber\\
\big(q^{-1}-q\big)\pi^+_k(F)  e_n  = -q^{\frac12-k-n}
\sqrt{(1-q^{2n})(1-q^{4k+2n-2})}  e_{n-1},\nonumber\\
\big(q^{-1}-q\big)^2 \pi^+_k(\Omega)  e_n  = \big(q^{2k-1}+q^{1-2k}-2\big)  e_n.\label{eq:pos}
\end{gather}

\textbf{Negative discrete series.} The negative discrete series
$\pi^-_k$ are labeled by $k>0$. The representation space is
$\ell^2(\N)$ with orthonormal basis $\{e_n\}_{n \in \N}$. The action is given by
\begin{gather}
\pi^-_k(K) e_n = q^{-(k+n)} e_n, \qquad \pi^-_k\big(K^{-1}\big) e_n =q^{k+n} e_n,
 \nonumber\\
 \big(q^{-1}-q\big)\pi^-_k(E)  e_n  = -q^{\frac12-k-n}
\sqrt{(1-q^{2n})(1-q^{4k+2n-2})}  e_{n-1},
\nonumber\\
\big(q^{-1}-q\big)\pi^-_k(F) e_n = q^{-\frac12-k-n}
\sqrt{(1-q^{2n+2})(1-q^{4k+2n})} e_{n+1}, \nonumber\\
\big(q^{-1}-q\big)^2 \pi^-_k(\Omega) e_n = \big(q^{2k-1}+q^{1-2k}-2\big) e_n.\label{eq:neg}
\end{gather}

\textbf{Principal unitary series.} The principal unitary series
representations $\pi^{\mathrm P}_{\rho,\eps}$ are labeled by $0\leq\rho <
-\frac{\pi}{2\ln q}$ and $\eps \in [0,1)$, where $(\rho,\eps)
\neq (0,\frac12)$. The representation space is $\ell^2(\Z)$ with
orthonormal basis $\{e_n\}_{n \in \Z}$. The action is given by
\begin{gather}
\pi^{\mathrm P}_{\rho,\eps}(K) e_n = q^{n+\eps} e_n, \qquad
\pi^{\mathrm P}_{\rho,\eps} \big(K^{-1}\big) e_n = q^{-(n+\eps)}e_n,\nonumber\\
\big(q^{-1}-q\big)\pi^{\mathrm P}_{\rho,\eps}(E) e_n = q^{-\frac12-n-\eps} \sqrt{(1-q^{2n+2\eps+2i\rho+1})
(1-q^{2n+2\eps-2i\rho+1})} e_{n+1},
\nonumber\\
\big(q^{-1}-q\big)\pi^{\mathrm P}_{\rho,\eps}(F) e_n = -q^{\frac12-n-\eps}
\sqrt{(1-q^{2n+2\eps+2i\rho-1}) (1-q^{2n+2\eps-2i\rho-1})} e_{n-1},
\nonumber\\
\big(q^{-1}-q\big)^2 \pi^{\mathrm P}_{\rho,\eps} (\Omega) e_n
=\big(q^{2i\rho}+q^{-2i\rho} -2\big) e_n.\label{eq:princ}
\end{gather}
For $(\rho,\eps)= (0,\frac12)$ the representation $\pi^{\mathrm P}_{0,\frac12}$ splits into a direct
sum of a positive and a negative discrete series representation: $\pi^{\mathrm P}_{0,\frac12}= \pi^+_\frac12
\oplus \pi^-_\frac12$. The representation space splits into two invariant subspaces: $\{ e_n\, | \, n \geq 0 \} \oplus \{e_n \, |\,
n<0 \}$.

\textbf{Complementary series.} The complementary series
representations $\pi^{\mathrm C}_{\lambda,\eps}$ are labeled by $\lambda$ and
$\eps$, where $\eps \in [0,\frac12)$ and $\lambda \in (-\frac12,-\eps)$, or
$\eps \in (\frac12,1)$ and $\lambda \in (-\frac12, \eps-1)$. The
representation space is $\ell^2(\Z)$ with orthonormal basis
$\{e_n\}_{n \in \Z}$. The action of the generators is the same as for the principal unitary series, with $i\rho$ replaced by $\lambda+\frac12$:
\begin{gather}
\pi^{\mathrm C}_{\lambda,\eps}(K) e_n = q^{n+\eps}\, e_n, \qquad
\pi^{\mathrm C}_{\lambda,\eps} \big(K^{-1}\big) e_n = q^{-(n+\eps)}e_n,\nonumber\\
\big(q^{-1}-q\big)\pi^{\mathrm C}_{\rho,\eps}(E) e_n
= q^{-\frac12-n-\eps} \sqrt{(1-q^{2n+2\eps+2\lambda+2})
(1-q^{2n+2\eps-2\lambda})} e_{n+1},
\nonumber\\
\big(q^{-1}-q\big)\pi^{\mathrm C}_{\lambda,\eps}(F) e_n = -q^{\frac12-n-\eps}
\sqrt{(1-q^{2n+2\eps+2\lambda}) (1-q^{2n+2\eps-2\lambda-2})} e_{n-1},\nonumber\\
\big(q^{-1}-q\big)^2 \pi^{\mathrm C}_{\lambda,\eps} (\Omega) e_n
=\big(q^{2\lambda+1}+q^{-1-2\lambda} -2\big) e_n.\label{eq:comp}
\end{gather}

\textbf{Strange series.} The f\/ifth class consists of the strange series representations $\pi^{\mathrm S}_{\sigma,\eps}$, labeled by $\sigma>0$ and $\eps \in [0,1)$, with representation space $\ell^2(\Z)$. The actions of the generators are the same as in principal unitary series with $i\rho$ replaced by $- \frac{ i\pi}{2\ln q}+\sigma$:
\begin{gather}
\pi^{\mathrm S}_{\sigma,\eps}(K) e_n = q^{n+\eps} e_n, \qquad
\pi^{\mathrm S}_{\sigma,\eps} \big(K^{-1}\big) e_n = q^{-(n+\eps)} e_n,\nonumber\\
\big(q^{-1}-q\big)\pi^{\mathrm S}_{\sigma,\eps}(E) e_n
= q^{-\frac12-n-\eps} \sqrt{(1+q^{2n+2\eps+2\sigma+1})
(1+q^{2n+2\eps-2\sigma+1})} e_{n+1},\nonumber\\
\big(q^{-1}-q\big)\pi^{\mathrm S}_{\sigma,\eps}(F) e_n = -q^{\frac12-n-\eps}
\sqrt{(1+q^{2n+2\eps+2\sigma-1}) (1+q^{2n+2\eps-2\sigma-1})} e_{n-1},\nonumber\\
\big(q^{-1}-q\big)^2 \pi^{\mathrm S}_{\sigma,\eps} (\Omega) e_n =-\big(q^{2\sigma}+q^{-2\sigma} +2\big)e_n.\label{eq:strange}
\end{gather}

The representations $\pi^{\mathrm P}_{\rho,\eps}$, $\pi^{\mathrm C}_{\lambda,\eps}$, $\pi^{\mathrm S}_{\sigma,\eps}$ with $\eps \in \R\setminus [0,1)$ are equivalent to the representations $\pi^{\mathrm P}_{\rho,\eps+m}$,  $\pi^{\mathrm C}_{\lambda,\eps+m}$, $\pi^{\mathrm S}_{\sigma,\eps+m}$ respectively, where $m \in \Z$ is such that $\eps+m \in [0,1)$. Therefore we will use these representations with any $\eps \in \R$. Furthermore, the principal unitary series representation~$\pi_{\rho,\eps}^{\mathrm P}$ is equivalent to $\pi_{-\rho,\eps}^{\mathrm P}$ and to $\pi_{\rho',\eps}^{\mathrm P}$ with $\rho'=\rho+\frac{m\pi}{\ln q}$, $m \in \Z$, and the complementary series representation $\pi_{\lambda,\eps}^{\mathrm{C}}$ is equivalent to $\pi_{-\lambda-1,\eps}^{\mathrm{C}}$.

It will be convenient to use the convention $e_{-n}=0$ for $n=1,2,\ldots$, if $e_n$ is a basis vector of $\ell^2(\N)$. Furthermore, we sometimes denote a basis vector of the representation space of the positive discrete series by $e_n^+$, and similarly for basis vectors of other representation spaces.

\section{A quantum analog of the tensor product of negative\\ and positive discrete series representations} \label{sec:-+}

Let $k_1,k_2>0$ and $\eps \in \mathbb R$. We consider the $\U_q$-representation
\begin{gather} \label{eq:defT1}
\mathcal T = \mathcal T_{k_1,k_2,\eps}= \Big(\big( \pi_{k_1}^- \otimes \pi_{k_2}^+\big) \oplus \big(\pi^{\mathrm S}_{k_1-\frac12,-\eps-k_1}  \otimes \pi_{k_2-\frac12,\eps+k_2}^{\mathrm S}\big)\Big)  \Delta.
\end{gather}
For the decomposition of $\mathcal T$ into irreducible representations we need the big $q$-Jacobi functions~\cite{KS03}. Let $a$, $b$, $c$ be parameters satisfying the conditions $a,b,c>0$ and $ab,ac,bc>1$. The big $q$-Jacobi functions are def\/ined by
\[
\Phi_\gamma(y;a,b,c|q)= \frac{ (\gamma/a, 1/bc, -y/abc\gamma;q)_\infty} { (1/ab, 1/ac, -y/bc;q)_\infty} \rphis{3}{2}{ 1/b\gamma, 1/c\gamma, -y/bc}{1/bc, -y/abc\gamma}{q, \gamma/a},
\]
for $y \in \C\setminus(-\infty, -bc]$ and $|\gamma|< a$. By \cite[(III.9)]{GRa04} $\Phi_\gamma$ is symmetric in $\gamma$, $\gamma^{-1}$, so for $|\gamma|\geq a$ we def\/ine the big $q$-Jacobi function by the same formula with $\gamma$ replaced by $\gamma^{-1}$. If $|y|<bc$, the $_3\varphi_2$-function can be transformed by \cite[(III.10)]{GRa04} and then we have
\begin{gather} \label{eq:3phi2}
\Phi_\gamma(y;a,b,c|q)= \rphis{3}{2}{\gamma/a, 1/a\gamma, -1/y}{1/ab, 1/ac}{q,-y/bc}.
\end{gather}
Let $t>0$. We def\/ine discrete sets $\Gamma^{\mathrm{f\/in}}$ and $\Gamma^{\mathrm{inf}}$ (`f\/in' and `inf' stand for `f\/inite' and `inf\/inite') by
\begin{gather*}
\Gamma^{\mathrm{f\/in}}=  \left\{ \frac{q^k}{e}\, \Big| \, e \in \{a,b,c\},\, k \in \N \text{ such that } \frac{q^k}{e}>1 \right\},\\
\Gamma^{\mathrm{inf}} = \left\{ -\frac{ abcq^k }{t}\, \Big| \, k \in \Z \text{ such that } \frac{abcq^k}{t}>1 \right\},
\end{gather*}
Note that the set $\Gamma^{\mathrm{f\/in}}$ is empty if $a,b,c>1$.
Now we introduce the measure $\nu(\,\cdot\,)=\nu(\,\cdot\,;a,b,c;t|q)$ by
\[
\int f(\gamma) d\nu(\gamma) = \frac1{4\pi i}\int_\T f(\gamma) v(\gamma) \frac{ d\gamma}{\gamma} + \sum_{\gamma \in \Gamma^{\mathrm{f\/in}} \cup \Gamma^{\mathrm{inf}}} f(\gamma)w(\gamma) ,
\]
where $v$ is the weight function on the counter-clockwise oriented unit circle $\T$ given by
\begin{gather*}
v(\gamma)=v(\gamma;a,b,c;t|q) =  \frac{ \theta(-t;q) }{\theta(-t/ab, -t/ac, -t/bc;q) (q,1/ab,1/ac,1/bc;q)_\infty}\nonumber\\
\phantom{v(\gamma)=v(\gamma;a,b,c;t|q) =}{}  \times \frac{ (\gamma^{\pm 2};q)_\infty }{ (\gamma^{\pm 1}/a, \gamma^{\pm 1}/b, \gamma^{\pm 1}/c;q)_\infty\, \theta(-t\gamma^{\pm 1}/abc;q)},
\end{gather*}
and{\samepage
\[
w(\gamma)=w(\gamma;a,b,c;t;q) = \Res{\gamma'= \gamma} \frac{ v(\gamma') }{\gamma'}, \qquad \gamma \in \Gamma^{\mathrm{f\/in}} \cup \Gamma^{\mathrm{inf}}.
\]
The weights $w$ in the discrete mass points can be calculated explicitly, see \cite[Section~8]{KS03}.}

Let $\mathcal H=\mathcal H(a,b,c;t|q)$ be the Hilbert space consisting of functions that satisfy $f(\gamma)=f(\gamma^{-1})$ ($\nu$-a.e.) with inner product
\[
\langle f ,g \rangle = \int f(\gamma) \overline{g(\gamma)}\, d\nu(\gamma).
\]
The set of functions $\{\gamma \mapsto \Phi_\gamma(x;a,b,c)\, | \, x \in \{-q^k\mid k \in \N\}\cup\{tq^k\mid k \in \Z \}\}$ is an orthogonal basis for $\mathcal H$, i.e.
\[
\langle \Phi_\cdot(y), \Phi_\cdot(y') \rangle = N_y \delta_{yy'}, \qquad y,y' \in \big\{{-}q^k\mid k \in \N\big\}\cup\big\{tq^k \mid k \in \Z\big\},
\]
where the quadratic norm $N_y$ is given by
\[
N_y =N_y(a,b,c;q)= \begin{cases}
 \dfrac{ (q,1/bc;q)_k }{(1/ab,1/ac;q)_k}q^{-k}, & y=-q^k, k \in \N ,\vspace{2mm}\\
 \dfrac{ (q,1/bc,-tq^k/ab, -tq^k/ac;q)_\infty }{ (1/ab,1/ac,-tq^k/bc, -tq^{k+1};q)_\infty }\frac{q^{-k}}{t}, & y=tq^k, k \in \Z.
\end{cases}
\]
We def\/ine
\[
\phi_y(x)=\phi_y(x;a,b,c|q) = \frac{\Phi_x(y;a,b,c|q)}{\sqrt{N_y(a,b,c;q)}},
\]
then $\{ \phi_{-q^k} \}_{k \in \N} \cup \{\phi_{tq^k} \}_{k \in \Z}$ is an orthonormal basis for $\mathcal H$. For ease of notations, it is useful to def\/ine $\phi_{-q^k}=0$ for $k=-1,-2,\ldots$.

The dif\/ference equation for the big $q$-Jacobi functions leads to a recurrence relation for the functions $\phi_{zq^k}$, $z\in\{-1,t\}$, of the form
\begin{gather} \label{eq:3termrecphi}
\big(x+x^{-1}\big)\phi_{zq^k}(x) = a_k \phi_{zq^{k+1}}(x) + b_k \phi_{zq^k}(x) + a_{k-1}\phi_{zq^{k-1}}(x),
\end{gather}
where
\begin{gather*}
a_k =a_k(a,b,c;z|q)=  \frac{ abcq^{-\frac12-2k} }{z^2} \sqrt{ \big(1+ zq^{k+1} \big) \left(1+\frac{ zq^k}{ab} \right) \left(1+\frac{zq^k}{ac} \right) \left(1+\frac{zq^k}{bc} \right) },\\
b_k =b_k(a,b,c;z|q)=  a+a^{-1} -\frac{abcq^{-2k}}{z^2}\left( 1+ \frac{zq^k}{ab} \right) \left( 1+ \frac{zq^k}{ac} \right)\\
\phantom{b_k =b_k(a,b,c;z|q)=}{}  -\frac{abcq^{1-2k}}{z^2} \big( 1+zq^k \big) \left( 1+\frac{zq^{k-1}}{bc} \right).
\end{gather*}
Here we assume $k \in \N$ if $z=-1$ and $k \in \Z$ if $z=t$, and we use $a_{-1}(a,b,c;-1|q)=0$. Note that $a_k$ and $b_k$ are both symmetric in $a,b,c$.

Let us remark that for $z=-1$ the recurrence relation corresponds to the three term recurrence relation for the continuous dual $q^{-1}$-Hahn polynomials, i.e.~Askey--Wilson polynomials with one parameter equal to zero and base $q^{-1}$;
\[
p_k(x)=p_k\big(x;a,b,c|q^{-1}\big) = a^{-k} \sqrt{ \frac{(ab,ac;q^{-1})_k}{(q^{-1}, bc;q^{-1})_k} } \rphis{3}{2}{q^{k}, ax, a/x}{ab, ac}{q^{-1}, q^{-1}}.
\]
So the function $\phi_{-q^k}$ is a multiple of the continuous dual $q^{-1}$-Hahn polynomials $p_k$. The moment problem corresponding to the continuous dual $q^{-1}$-Hahn polynomials is indeterminate, so the measure $\nu$ def\/ined above is a (non-extremal) solution for this moment problem.

We now turn to the decomposition of $\mathcal T$. First we diagonalize $\mathcal T(\Omega)$ on a suitable subspace of $\ell^2(\N)^{\otimes 2} \oplus \ell^2(\Z)^{\otimes 2}$. We def\/ine subspaces $\mathcal A_p$, $p\in\Z$, by
\[
\mathcal A_p =  \mathrm{span}\{f_{p,n}^-\ | \ n\in \N\}  \oplus  \mathrm{span} \{f_{p,n}^+\ | \ n \in \Z \},
\]
where
\[
f_{p,n}^- = \begin{cases}
e_n \otimes e_{n-p}, & p \leq 0,\\
e_{n+p} \otimes e_n, & p \geq 0,
\end{cases}
 \qquad
f_{p,n}^+=
\begin{cases}
(-1)^n e_{-n} \otimes e_{p+n}, & p \leq 0,\\
(-1)^{n+p} e_{-n-p} \otimes e_n, & p \geq 0.
\end{cases}
\]
It is useful to observe that $f_{p,n}^-=0$ for $n=-1,-2,\ldots$. Furthermore, note that $\overline{\mathcal A_p} \cong \ell^2(\N) \oplus \ell^2(\Z)$ and $\mathcal A = \bigoplus_{p \in \Z} \mathcal A_p$ is a dense subspace of $\ell^2(\N)^{\otimes 2} \oplus \ell^2(\Z)^{\otimes 2}$. From \eqref{eq:De(OM)} it follows that each subspace $\mathcal A_p$ is invariant under $\mathcal T(\Omega)$.

\begin{proposition} \label{prop:Te}
Let $p \in \Z$, $\sigma \in \{-,+\}$, and let $a$, $b$, $c$, $t_\sigma$ be defined by
\begin{gather}
a= \begin{cases}
q^{2k_2-2k_1-2p-1}, & p \geq 0,\\
q^{2k_1-2k_2+2p-1}, & p \leq 0,
\end{cases}
 \qquad b = q^{1-2k_1-2k_2}, \qquad
c=\begin{cases}
q^{2k_1-2k_2-1}, & p \geq 0,\\
q^{2k_2-2k_1-1}, & p \leq 0,
\end{cases}
\nonumber\\  t_\sigma =
\begin{cases}
-1, & \sigma=-,\\
q^{2\eps}, & \sigma=+.
\end{cases}\label{eq:abct}
\end{gather}
Then the operator $\Theta:\mathcal A_p \rightarrow \mathcal H(a,b,c;t_+|q^2)$ defined by
\[
f_{p,n}^\sigma \mapsto \left(x\mapsto \phi_{t_\sigma q^{2n}}\big(x;a,b,c|q^2\big) \right), \qquad n \in \Z,
\]
extends to a unitary operator. Moreover, $\Theta$ intertwines $\mathcal T(\Omega)$ with the multiplication operator $(q^{-1}-q)^{-2}M_{x+x^{-1}-2}$.
\end{proposition}

\begin{proof}
First of all, $\Theta$ extends to a unitary operator since it maps one orthonormal basis onto another.

To check the intertwining property, we consider the explicit action of $\Omega$.
First assume $p \geq 0$. From \eqref{eq:De(OM)}, \eqref{eq:pos}, \eqref{eq:neg} and \eqref{eq:strange} it follows that the action of the Casimir operator is given by
\[
\big(q^{-1}-q\big)^2 [\mathcal T(\Omega)+2 ] f_{p,n}^\sigma = a_n^\sigma f_{p,n+1}^\sigma + b_n^\sigma f_{p,n}^\sigma + a_{n-1}^\sigma f_{p,n-1}^\sigma,
\]
where
\begin{gather*}
a_n^\sigma =   t_\sigma^{-2}q^{-2-2k_2-2k_1-2p-4n}\\
 \hphantom{a_n^\sigma =}{}
 \times \sqrt{ \big(1+t_\sigma q^{4k_1+2p+2n}\big)\big(1+t_\sigma q^{2p+2n+2}\big)\big(1+t_\sigma q^{4k_2+2n}\big) \big(1+t_\sigma q^{2n+2}\big)},\\
b_n^\sigma=   q^{-2p-2k_1+2k_2-1} +q^{2p+2k_1-2k_2+1} - t_\sigma^{-2}q^{1-2k_1-2k_2-2p-4n}\big(1+t_\sigma q^{2n+4k_2-2}\big)\big(1+t_\sigma q^{2n}\big) \\
\phantom{b_n^\sigma=}{}
 + t_\sigma^{-2}q^{-1-2k_1-2k_2-2p-4n}\big(1+t_\sigma q^{4k_1+2p+2n}\big)\big(1+t_\sigma q^{2p+2n+2}\big).
\end{gather*}
The intertwining property follows from comparing this with the recurrence relation \eqref{eq:3termrecphi} for the function $\phi_{t_\sigma q^{2n}}(x;a,b,c|q^2)$. For $p <0$ the proof runs along the same lines.
\end{proof}

By Proposition~\ref{prop:Te} $\Theta$ intertwines $\mathcal T(\Omega)$ with a multiplication operator on~$\mathcal H$. This allows us to read of\/f the spectrum of $\mathcal T(\Omega)$ from the support of the measure $\nu$. Since the irreducible $*$-representations are characterized by the actions of $\Omega$ and $K$, we now only need to consider the action of $K$ to f\/ind the decomposition of $\mathcal T$ into irreducible representations;
\[
\mathcal T(K) f_{p,n}^\sigma = q^{-2p + 2k_2 -2k_1}  f_{p,n}^\sigma.
\]
By comparing the spectral values of $\mathcal T(K)$ and $\mathcal T(\Omega)$ with \eqref{eq:pos}--\eqref{eq:strange}, we obtain the following decomposition of $\mathcal T$.

\begin{theorem} \label{thm:decompT1}
The $\U_q$-representation $\mathcal T_{k_1,k_2,\eps}$ is unitarily equivalent to:
\begin{alignat*}{4}
& (i)\quad && \dirint{0}{-\frac{\pi}{2 \ln q}} \pi^{\mathrm P}_{\rho,\eps'} d\rho \oplus \bigoplus_{\substack{j \in \Z \\ \sigma_j>0}} \pi^{\mathrm S}_{\sigma_j,\eps'},\qquad  &&  |k_1-k_2|\leq \frac12,\ k_1+k_2\geq \frac12,& \\
& (ii)\quad && \dirint{0}{-\frac{\pi}{2 \ln q}} \pi^{\mathrm P}_{\rho,\eps'} d\rho \oplus \bigoplus_{\substack{j \in \Z \\ \sigma_j>0}} \pi^{\mathrm S}_{\sigma_j,\eps'} \oplus \pi^{\mathrm C}_{-k_1-k_2,\eps'}, \qquad && k_1+k_2<\frac12,& \\
& (iii)\quad && \dirint{0}{-\frac{\pi}{2 \ln q}} \pi^{\mathrm P}_{\rho,\eps'} d\rho \oplus \bigoplus_{\substack{j \in \Z \\ \sigma_j>0}} \pi^{\mathrm S}_{\sigma_j,\eps'} \oplus \bigoplus_{\substack{j \in \N\\k_j^+>\frac12}} \pi^+_{k_j^+} , \qquad && k_1-k_2+\frac12<0,& \\
& (iv)\quad && \dirint{0}{-\frac{\pi}{2 \ln q}} \pi^{\mathrm P}_{\rho,\eps'} d\rho \oplus
\bigoplus_{\substack{j \in \Z \\ \sigma_j>0}} \pi^{\mathrm S}_{\sigma_j,\eps'} \oplus \bigoplus_{\substack{j \in \N\\k_j^->\frac12}} \pi^-_{k_j^-} , \qquad && k_2-k_1+\frac12<0,&
\end{alignat*}
where $\eps'=k_2-k_1$, $\sigma_j= k_1+k_2+\eps+\frac12+j$, $k_j^+=k_2-k_1-j$ and $k_j^-=k_1-k_2-j$. For $n,p \in \Z$ and $\sigma \in \{-,+\}$ the unitary intertwiner is given by
\begin{alignat*}{3}
& (i)\quad && f_{p,n}^\sigma \mapsto \Lambda f_{p,n}^\sigma,  & \\
& (ii) \quad && f_{p,n}^\sigma \mapsto \Lambda f_{p,n}^\sigma + (\Theta f_{p,n}^\sigma)\big(q^{2k_1+2k_2-1}\big)\sqrt{w\big(q^{2k_1+2k_2-1};a,b,c;t_+;q^2\big)} e_{-p}^{\mathrm C},  & \\
& (iii) \quad && f_{p,n}^\sigma \mapsto \Lambda f_{p,n}^\sigma +
\begin{cases}
\displaystyle \sum_{\substack{j \in \N\\k_{j+p}^+>\frac12}} (\Theta f_{p,n}^\sigma)\big(q^{2k_{j+p}^+-1}\big) \sqrt{w\big(q^{2k_{j+p}^+-1};a,b,c;t_+;q^2\big)} e_{j}^+, & p \geq 0,\\
\displaystyle \sum_{\substack{j \in \N\\k_{j}^+>\frac12}} (\Theta f_{p,n}^\sigma)\big(q^{2k_{j}^+-1}\big) \sqrt{w\big(q^{2k_{j}^+-1};a,b,c;t_+;q^2\big)} e_{j-p}^+, & p \leq 0,
\end{cases} &
\\
& (iv) \quad && f_{p,n}^\sigma \mapsto \Lambda f_{p,n}^\sigma +
\begin{cases}
\displaystyle \sum_{\substack{j \in \N\\k_j^->\frac12}} (\Theta f_{p,n}^\sigma)\big(q^{2k_j^--1}\big) \sqrt{w\big(q^{2k_j^--1};a,b,c;t_+;q^2\big)}e_{j+p}^-, & p \geq 0,\\
\displaystyle \sum_{\substack{j \in \N\\k_{j-p}^->\frac12}} (\Theta f_{p,n}^\sigma)\big(q^{2k_{j-p}^--1}\big) \sqrt{w\big(q^{2k_j^--1};a,b,c;t_+;q^2\big)} e_{j}^-, & p \leq 0,
\end{cases}&
\end{alignat*}
where
\begin{gather*}
\Lambda f_{p,n}^\sigma =  \int_0^{-\frac{\pi}{2\ln q}} \big(\Theta f_{p,n}^\sigma\big)\big(q^{2i\rho}\big) \sqrt{v\big(q^{2i\rho};a,b,c;t_+|q^2\big)}  e_{-p}^{\mathrm P}\, d\rho\\
\phantom{\Lambda f_{p,n}^\sigma =}{} + \sum_{\substack{j \in \Z\\\sigma_j>0}} \big(\Theta f_{p,n}^\sigma\big)\big({-}q^{2\sigma_j}\big) \sqrt{w\big({-}q^{2\sigma_j};a,b,c;t_+;q^2\big)} e_{-p}^{\mathrm S},
\end{gather*}
$\Theta$ is as in Proposition {\rm \ref{prop:Te}} and $a$, $b$, $c$, $t_+$ are given by \eqref{eq:abct}.
\end{theorem}

Note that in case $(i)$ the intertwiner maps from $\ell^2(\N)^{\otimes 2} \oplus \ell^2(\Z)^{\otimes 2}$ into $\int_0^{-\pi/(2\ln(q))} \ell^2(\Z)\, d\rho \oplus \bigoplus_{\N} \ell^2(\Z)$. In case $(ii)$ another $\ell^2(\Z)$ has to be added here, and in cases $(iii)$ and $(iv)$ a f\/inite number of $\ell^2(\N)$-spaces has to be added.

\begin{remark}
$(i)$ Letting $q \uparrow 1$ we (formally) obtain the decomposition of the tensor product of a negative and a positive discrete series representation of $\su(1,1)$, see e.g.~\cite[Theorem 7.3]{Re} or \cite[Theorem 2.2]{GK02}.

$(ii)$ Theorem \ref{thm:decompT1} shows that the big $q$-Jacobi functions have an interpretation as quantum analogs of Clebsch--Gordan coef\/f\/icients. In the classical case the Clebsch--Gordan coef\/f\/icients for the tensor product of negative and positive discrete series are essentially continuous dual Hahn polynomials, see \cite[Section~4]{OZ} or \cite[Theorem~2.2]{GK02}. So in the context of Clebsch--Gordan coef\/f\/icients the big $q$-Jacobi transform pair should be considered as a $q$-analog of the transform pair corresponding to continuous dual Hahn polynomials.

$(iii)$ The $\U_q$-representation $\pi_{k_2}^+\otimes \pi_{k_1}^-$ can be decomposed similar as in Theorem \ref{thm:decompT1}, but the inf\/inite sum of strange series representations does not occur in this situation. The Clebsch--Gordan coef\/f\/icients in this case are essentially continuous dual $q^2$-Hahn polynomials, see \cite[Section 2]{KM} or \cite[Theorem 2.4]{Gr04}.

$(iv)$ The term $\big( \pi_{k_1}^- \otimes \pi_{k_2}^+\big)\Delta$ occurring in the def\/inition of $\mathcal T$ would of course also be a~quantum analogue of the tensor product of a negative and positive discrete series representation of~$\su(1,1)$, but we consider it unlikely that this representation can be decomposed into irreducible representations. Indeed, in this case the action of the Casimir operator corresponds to the Jacobi operator for the continuous dual $q^{-2}$-Hahn polynomials, see also \cite[Remark~8.1]{Gr06}. The corresponding moment problem is indeterminate (so the Casimir operator is not essentially self-adjoint in this case!) and no explicit $N$-extremal solutions are known. Even if such a measure was known, it would have discrete support, implying that the decomposition would be a~direct sum of irreducible representations, and not a direct integral as in the classical case.
\end{remark}

Let us denote the intertwiner from Theorem \ref{thm:decompT1} by $I$. The actions $I\circ \mathcal T(X)\circ I^{-1}$, $X=E,F$, are given by the appropriate actions of $E$ and $F$ in \eqref{eq:pos}--\eqref{eq:strange}, up to a phase factor. It is possible to determine the actions explicitly using the explicit expressions for the weight functions~$v$ and~$w$, the explicit expressions for $\Theta f_{p,n}^\sigma$ as $_3\varphi_2$-functions, and the following contiguous relations for $_3\varphi_2$-functions, see \cite[(2.3), (2.4)]{GIM96},
\begin{gather*}
\rphis{3}{2}{Aq,B,C}{D,E}{q,\frac{DE}{ABCq}} - \rphis{3}{2}{A,B,C}{D,E}{q,\frac{DE}{ABC}}  \\
\qquad{}=\frac{DE(1-B)(1-C)}{ABCq(1-D)(1-E)} \rphis{3}{2}{Aq,Bq,Cq}{Dq,Eq}{q,\frac{DE}{ABCq}},
\\
\left(1-\frac{D}{A}\right)\left(1-\frac{E}{A}\right) \rphis{3}{2}{A/q,B,C}{D,E}{q,\frac{DEq}{ABC}} \\
\qquad{} -\left(1-\frac{q}{A}\right)\left(1-\frac{DE}{ABC}\right)\rphis{3}{2}{A,B,C}{D,E}{q,\frac{DE}{ABC}}
\\
\qquad{}=\frac{q}{A}\left(1-\frac{D}{q} \right) \left(1-\frac{E}{q} \right)\rphis{3}{2}{A/q,B/q,C/q}{D/q,E/q}{q,\frac{DEq}{ABC}}.
\end{gather*}
We do not work out the details.

\subsection[Completing $\mathcal T$ to a genuine tensor product representation]{Completing $\boldsymbol{\mathcal T}$ to a genuine tensor product representation}

In this subsection we def\/ine a representation $\mathcal T'$ that completes $\mathcal T$ def\/ined by \eqref{eq:defT1} to a genuine tensor product. Let $k_1,k_2>0$ and $\eps \in \frac12\Z$, and def\/ine the $\U_q$-representation $\mathcal T'$ by
\begin{gather*} 
\mathcal T' = \mathcal T_{k_1,k_2,\eps}'= \Big( \big(\pi^-_{k_1}\otimes \pi^{\mathrm S}_{k_2-\frac12,\eps+k_2} \big)\oplus \big( \pi^{\mathrm S}_{k_1-\frac12,-\eps-k_1} \otimes \pi^+_{k_2}\big)\Big)  \Delta.
\end{gather*}
Now the sum $\mathcal T \oplus \mathcal T'=\Big( \big(\pi^-_{k_1}\oplus \pi^{\mathrm S}_{k_1-\frac12,-\eps-k_1}\big) \otimes \big( \pi^+_{k_2} \oplus\pi^{\mathrm S}_{k_2-\frac12,\eps+k_2} \big)\Big)\, \Delta$ is a genuine tensor product representation of $\U_q$, which can also be considered as a quantum analog of the tensor product of a negative and a positive discrete series representation of~$\su(1,1)$.

The decomposition of $\mathcal T'$ into irreducible $*$-representations is established in the same way as the decomposition of $\mathcal T$, therefore we omit most of the details. We remark that we need here the condition $\eps \in \frac12\Z$ (instead of $\eps \in \R$), because our method for constructing an intertwiner only works if basis vectors of $\ell^2(\Z)$ can be labeled by $2\eps+m$ for $m \in \Z$, which forces $\eps$ to be in~$\frac12\Z$.

For the diagonalization of $\mathcal T'(\Omega)$ we need the big $q$-Jacobi polynomials \cite{AA}, \cite[\S~14.5]{KLS}. They are def\/ined by
\[
P_m(x;a,b,c;q) = \rphis{3}{2}{q^{-m},abq^{m+1},x}{aq,cq}{q,q}.
\]
We assume that $m \in \N$, $x \in \{aq^{k+1} \mid k \in \N\} \cup \{ cq^{k+1} \mid k \in \N\}$, $0<a,b<q^{-1}$ and $c<0$, then the big $q$-Jacobi polynomials satisfy the orthogonality relations
\[
\int_{cq}^{aq} P_{m_1}(x)P_{m_2}(x) u(x) d_q x = \frac{\delta_{m_1m_2}}{v(m_1)}, \qquad \sum_{m=0}^\infty P_m(x_1)P_m(x_2) v(m) = \frac{\delta_{x_1x_1}}{|x_1|u(x_1)},
\]
where the (positive) functions $u(\cdot) = u(\cdot;a,b,c;q)$ and  $v(\cdot)=v(\cdot;a,b,c;q)$ are given by
\begin{gather*}
u(x;a,b,c;q) = \frac{ (x/a,x/c;q)_\infty }{(x,bx/c;q)_\infty},\nonumber \\
v(m;a,b,c;q) = \frac{ (1-abq^{2m+1})}{aq(1-abq)} \frac{ (aq,bq,cq,abq/c;q)_\infty }{ (q,abq^2;q)_\infty \theta(c/a) }\nonumber\\
 \phantom{v(m;a,b,c;q) =}{}
 \times \frac{ (aq,abq,cq;q)_m }{ (q,bq,abq/c;q)_m } \big({-}acq^2\big)^{-m} q^{-\frac12 m(m-1)}.
\end{gather*}
Here the Jackson $q$-integral is def\/ined by
\[
\int_{cq}^{aq} f(x)d_q x= (1-q)\sum_{k=0}^\infty f\big(aq^{k+1}\big)aq^{k+1} - (1-q)\sum_{k=0}^\infty f\big(cq^{k+1}\big)cq^{k+1}.
\]
We def\/ine functions $r_x(m)$, related to the big $q$-Jacobi polynomials, by
\begin{gather*} 
r_x(m;a,b,c;q) = \sqrt{|x| u(x;a,b,c;q) v(m;a,b,c;q)} P_m(x;a,b,c;q).
\end{gather*}
By the orthogonality relations for the big $q$-Jacobi polynomials we have
\[
\sum_{n=0}^\infty r_{x_1}(m) r_{x_2}(m)  = \delta_{x_1x_2}, \qquad
\sum_{k=0}^\infty  r_{aq^k}(m_1) r_{aq^k}(m_2) + \sum_{k=0}^\infty  r_{cq^k}(m_1) r_{cq^k}(m_2) = \delta_{m_1m_2}.
\]
Furthermore, from the $q$-dif\/ference equation for $P_m$ it follows that the functions $r_x(m)$ satisfy the following $q$-dif\/ference equation in $x$:
\[
\frac1{\sqrt{abq}}\big(1-q^{-m}\big)\big(1-abq^{m+1}\big) r_x(m) = A(x) r_{qx}(m) + B(x) r_x(m) + A(x/q) r_{x/q}(m),
\]
where
\begin{gather*}
A(x)  = -\frac{c\sqrt a}{ x^2\sqrt{b}} \sqrt{ (1-x)(1-x/a)(1-x/c)(1-bx/c) },\\
B(x)  = \frac{c\sqrt{aq}}{x^2\sqrt{b}}\big( (1-x)(1-bx/c)+q(1-x/aq)(1-x/cq) \big) \\
\phantom{B(x)}{}  = \frac{c\sqrt{aq}}{x^2\sqrt{b}}\big( (1-x)(1-x/a)+q(1-bx/cq)(1-x/cq) \big) \\
\phantom{B(x)=}{}
+ \sqrt{abq} + \frac{1}{\sqrt{abq}}-c\sqrt{\frac{q}{ab}}- \frac1c\sqrt{\frac{ab}{q}}.
\end{gather*}

Now we are ready to diagonalize $\mathcal T'(\Omega)$. For $p \in \Z$ and $n \in \N$ we def\/ine
\[
f_{p,n}^- = e_{n} \otimes e_{n+p}, \qquad f_{p,n}^+ = e_{2\eps-n+p} \otimes e_n,
\]
and we write
\[
\mathcal A_p = \mathrm{span}\{f_{p,n}^- \mid n \in \N\} \oplus \mathrm{span}\{ f_{p,n}^+ \mid n \in \N\}.
\]
Note that $\bigoplus_{p \in \Z} \mathcal A_p$ is a dense subspace of $\ell^2(\N)\otimes \ell^2(\Z) \oplus \ell^2(\Z)\otimes \ell^2(\N)$.

\begin{proposition} \label{prop:Te2}
Let $p \in \Z$ and define $\Theta:\mathcal A_p \to \ell^2(\N)$ by
\begin{gather*}
f_{p,n}^-   \mapsto \big(m \mapsto r_{q^{2n+4k_1}}\big(m;q^{4k_1-2},q^{4k_2-2}, -q^{4k_1-2p-2\eps-2};q^2\big)\big),\\
f_{p,n}^+  \mapsto \big(m \mapsto (-1)^n r_{-q^{2n+4k_1-2p-2\eps}}\big(m;q^{4k_1-2},q^{4k_2-2}, -q^{4k_1-2p-2\eps-2};q^2\big)\big),
\end{gather*}
then $\Theta$ extends to a unitary operator. Furthermore, $\Theta$ intertwines $\mathcal T'(\Omega)$ with the multiplication operator $-(q^{-1}-q)^{-2} M_{q^{2k_1+2k_2-1+2m}+q^{-(2k_1+2k_2-1+2m)}+2}$.
\end{proposition}

\begin{theorem} \label{thm:decompT2}
The $\U_q$-representation $\mathcal T'_{k_1,k_2,\eps}$ is unitarily equivalent to
\[
\bigoplus_{m=0}^\infty \pi^{\mathrm S}_{k_1+k_2-\frac12+m,\eps-k_1+k_2}.
\]
The unitary intertwiner $\ell^2(\N)\otimes \ell^2(\Z) \oplus \ell^2(\Z)\otimes \ell^2(\N) \to \bigoplus_{m=0}^\infty \ell^2(\Z)$ is given by
\[
f_{p,n}^\sigma \mapsto \sum_{m=0}^\infty \big(\Theta f_{p,n}^\sigma\big)(m)  e_p^{\mathrm S}, \qquad p \in \Z, \ n \in \N, \ \sigma \in \{-,+\}.
\]
\end{theorem}

\begin{remark}
Theorem \ref{thm:decompT2} shows that the big $q$-Jacobi polynomials have a natural interpretation as Clebsch--Gordan coef\/f\/icients for $\U_q$-representations. In this interpretation they do not have a classical analog, since the $\U_q$-representation $\mathcal T'$ vanishes in the classical limit.
\end{remark}

\section{A quantum analog of the tensor product\\ of two principal unitary series representations} \label{sec:PP}

Let $\rho_1,\rho_2 \in (0,-\frac{\pi}{2\ln q})$ and $\eps_1,\eps_2,\eps\in \R$. In this section we consider the representation
\[
\mathcal T = \mathcal T_{\rho_1,\rho_2,\eps_1,\eps_2,\eps} = \Big((\pi_{\rho_1,\eps_1+\eps}^{\mathrm P} \otimes \pi_{\rho_2,\eps_2-\eps}^{\mathrm P}) \oplus  (\pi_{\rho_1',\eps_1-\eps}^{\mathrm P} \otimes \pi_{\rho_2',\eps_2+\eps}^{\mathrm P})\Big)\Delta,
\]
where $\rho_j'=\rho_j-\frac{\pi}{2\ln q}$, $j=1,2$. The representation space is $\ell^2(\Z)^{\otimes 2} \oplus \ell^2(\Z)^{\otimes 2}$. Observe that for $\eps_1,\eps_2,\eps \in \frac12\Z$ the representation $\mathcal T_{\rho_1,\rho_2',\eps_1,\eps_2,\eps} \oplus \mathcal T_{\rho_1,\rho_2,\eps_1,\eps_2,\eps}$ is a genuine tensor product representation, using the equivalence $\pi_{\rho,\eps}^{\mathrm P}\cong\pi_{\rho,\eps+1}^{\mathrm P}$.

For the decomposition of $\mathcal T$ into irreducible representations we need the vector-valued big $q$-Jacobi functions~\cite{Gr09}. Let $a,c \in \C$, $z_+>0$ and $z_-<0$. We set $b=\overline{a}$, $d=\overline{c}$ and $s= \sqrt{q}|\frac{c}{a}|$. For a function $f$ depending on the parameter $a$ and $b$, $f = f(\,\cdot\,;a,b)$, we write
\[
f^\dagger=f^\dagger(\,\cdot\,;a,b)= f(\,\cdot\,;b,a).
\]
We def\/ine a discrete set $\Gamma$ depending on the parameters $a$, $c$, $z_+$, $z_-$ by
\begin{gather*}
\Gamma  =  \Gamma^{\mathrm{f\/in}}_{s} \cup \Gamma^{\mathrm{f\/in}}_{q/s} \cup \Gamma^{\mathrm{inf}},\\
\Gamma^{\mathrm{f\/in}}_\alpha  =  \left\{ \frac{1}{\alpha q^{k} }\ \Big| \ k\in\N,\ \alpha q^{k}>1 \right\},\\
\Gamma^{\mathrm{inf}}  = \left\{ z_-z_+q^{k-\frac12}|ac| \ \Big| \ k \in \Z, \ -z_-z_+q^{k-\frac12}|ac|<1 \right\}.
\end{gather*}

Let $\mathbf v$ be the matrix-valued function given by
\begin{gather*}
\gamma \mapsto \mathbf v(\gamma) =
\begin{pmatrix}
v_2(\gamma) & v_1^\dagger(\gamma) \\
v_1(\gamma) & v_2(\gamma)
\end{pmatrix},
\\
v_1(\gamma) = \frac{ (cq/a, dq/a;q)_\infty^2 \theta(bz_+,bz_-;q) }{(1-q)abz_-^2z_+^2 \theta(z_-/z_+,z_+/z_-,a/b,b/a;q)}\\
\phantom{v_1(\gamma) =}{} \times \frac{(\gamma^{\pm2};q)_\infty}{(s\gamma^{\pm1}, cq\gamma^{\pm1}/as, dq\gamma^{\pm1}/as;q)_\infty \theta(s\gamma^{\pm1},absz_- z_+\gamma^{\pm 1};q)}\\
\phantom{v_1(\gamma) =}{} \times \Big( z_-\theta(az_+, cz_+, dz_+, bz_-, asz_-\gamma^{\pm1};q) -  z_+\theta(az_-, cz_-, dz_-, bz_+, asz_+\gamma^{\pm1};q)  \Big),\\
v_2(\gamma) = \frac{(cq/a,dq/a,cq/b,dq/b;q)_\infty \theta(az_+, az_-, bz_+, bz_-, cdz_- z_+;q) }{abz_-^2z_+(1-q)\theta(z_+/z_-,a/b,b/a;q)} \\
\phantom{v_2(\gamma) =}{} \times \frac{ (\gamma^{\pm2};q)_\infty }{ (s\gamma^{\pm1};q)_\infty \theta(s\gamma^{\pm1}, absz_- z_+\gamma^{\pm1};q) }.
\end{gather*}
For $\gamma \in \T$ we have $v_1^\dagger(\gamma)=\overline{v_1(\gamma)}$. In this case the matrix $\mathbf v(\gamma)$ is positive def\/inite and we can write
\[
\mathbf v(\gamma) = S(\gamma)^T S(\gamma),
\]
where $S^T$ is the conjugate transpose of $S$, and the matrix-valued function $S$ is given by
\[
S(\gamma)=S(\gamma;a,c|q) =
\begin{pmatrix}
\sqrt{\dfrac{v_2(\gamma)+|v_1(\gamma)|}{2}} \dfrac{ v_1(\gamma) }{|v_1(\gamma)|} & \sqrt{\dfrac{v_2(\gamma)+|v_1(\gamma)|}{2}} \vspace{2mm}\\
- \sqrt{\dfrac{v_2(\gamma)-|v_1(\gamma)|}{2}} \dfrac{ v_1(\gamma) }{|v_1(\gamma)|} & \sqrt{\dfrac{v_2(\gamma)-|v_1(\gamma)|}{2}}
\end{pmatrix}.
\]
We also def\/ine
\begin{gather*}
v(\gamma) = v(\gamma;a,c;z_-,z_+|q)=  -\frac{\theta(cz_-, dz_-, cz_+, dz_+;q)}{z_+(1-q) \theta(z_-/z_+) } \\
\phantom{v(\gamma) = }{} \times
 \frac{(\gamma^2,q\gamma^2;q)_\infty}{ (cq\gamma/as, dq\gamma/as, cq\gamma/bs, dq\gamma/bs, s\gamma, q\gamma/s;q)_\infty  \theta(absz_- z_+/q\gamma;q) }.
\end{gather*}
With this function we def\/ine a positive weight function $w$ on $\Gamma$ by
\[
w(\gamma)=w(\gamma;a,c;z_-,z_+|q) = \frac1{b(\gamma)}\Res{\gamma'=\gamma} \frac{v(\gamma')}{\gamma'},\qquad \gamma \in \Gamma,
\]
where
\[
b(\gamma)=
\begin{cases}
\left(\dfrac{z_+}{z_-}\right)^{k+1}\dfrac{\theta(cz_-,dz_-;q)}{\theta(cz_+,dz_+;q)}, &\gamma =z_-z_+q^{k-\frac12}|ac| \in \Gamma^{\mathrm{inf}},\\
\left(\dfrac{z_-}{z_+}\right)^k ,& \gamma = s^{-1}q^{-k} \in \Gamma^{\mathrm{f\/in}}_s,\\
\left( \dfrac{z_-}{z_+} \right)^k \dfrac{ \theta(az_+, bz_+, cz_-, dz_-;q) }{ \theta(az_-, bz_-, cz_+, dz_+;q) }, &\gamma = sq^{-1-k} \in \Gamma^{\mathrm{f\/in}}_{q/s}.
\end{cases}
\]
Explicit expressions for $w$ can be obtained by evaluating the residues, see \cite[Section~4]{Gr09}.
We denote by $F(\T\cup \Gamma)$ the vector space of functions that are $\C^2$-valued on $\T$ and $\C$-valued on $\Gamma$. We def\/ine $\mathcal H=\mathcal H(a,c;z_-,z_+|q)$ to be the Hilbert space consisting of functions in $F(\T\cup\Gamma)$ satisfying $f(\gamma)=f(\gamma^{-1})$ almost everywhere on $\T$, that have f\/inite norm with respect to the inner product
\[
\langle f,g \rangle = \frac{ 1}{4\pi i } \int_\T g(\gamma)^T
\mathbf v(\gamma) f(\gamma) \frac{ d\gamma}{\gamma} + \sum_{\gamma \in \Gamma}  f(\gamma) \overline{g(\gamma)} w(\gamma).
\]

The vector-valued big $q$-Jacobi functions are def\/ined by
\begin{gather*} 
\Psi(y,\gamma) = \Psi(y,\gamma;a,c|q)=
\begin{cases}
\displaystyle \vect{\varphi_\gamma(y)}{\varphi_\gamma^\dagger(y)}, & \gamma \in \T,\\
d(\gamma) \varphi_\gamma(y)+d^\dagger(\gamma) \varphi_\gamma^\dagger(y), & \gamma \in \Gamma^{\mathrm{inf}} \cup \Gamma^{\mathrm{f\/in}}_{q/s},\\
c(\gamma) \varphi_\gamma(y), & \gamma \in \Gamma_s^{\mathrm{f\/in}}.
\end{cases}
\end{gather*}
Here $\varphi_\gamma$ is the reparametrized big $q$-Jacobi function given by
\[
\varphi_\gamma(y) = \varphi_\gamma(y;a,c|q) = \Phi_\gamma(-1/by; 1/s, as/cq, as/dq|q),
\]
so explicitly $\varphi_\gamma$ is given by, see \eqref{eq:3phi2},
\[
\varphi_\gamma(y) = \rphis{3}{2}{q/ay, s\gamma, s/\gamma}{cq/a, dq/a}{q,by}, \qquad |y| <  \frac{1}{|b|},
\]
and furthermore
\begin{gather*}
d(\gamma) = \frac{  (cq/a, dq/a;q)_\infty \theta(bz_+) }{  \theta(a/b,cz_+,dz_+) }   \frac{ (cq\gamma/sb, dq\gamma/sb;q)_\infty \theta(asz_+/q\gamma) }{ (q\gamma^2,s/\gamma;q)_\infty},\\
c(\gamma) = \frac{ (cq/a, dq/a, 1/\gamma^2;q)_\infty \theta(bz) }{ (s/\gamma, cq/as\gamma, dq/as\gamma;q)_\infty \theta( bsz\gamma ) }.
\end{gather*}
For $\gamma \in \Gamma^{\mathrm{f\/in}}_s \cup \Gamma^{\mathrm{f\/in}}_{q/s}$ the function $\Psi(\,\cdot\,;\gamma)$ is actually a multiple of a big $q$-Jacobi polynomial, see \cite[Lemma 3.9]{Gr09}.

The set
\[
\left\{\gamma \mapsto \Psi(y,\gamma) \mid y \in \big\{z_-q^k \mid k \in \Z\big\} \cup \big\{z_+q^k \mid k \in \Z\big\} \right\}
\]
is an orthogonal basis for $\mathcal H$. We have
\[
\langle \Psi(y,\,\cdot\,),\Psi(y',\,\cdot\,)\rangle = N_y\delta_{yy'},
\]
where the squared norm $N_y$ is given by
\[
N_y=N_y(a,c;q) =\frac1{|y|}\left| \frac{(cy;q)_\infty}{(ay;q)_\infty } \right|^2.
\]
We def\/ine
\[
\psi_y(x)=\psi_y(x;a,c|q) = \frac{\Psi(y,x;a,c|q)}{\sqrt{N_y(a,c;q)}},
\]
then $\{\psi_{z_-q^k}\}_{k \in \Z} \cup \{\psi_{z_+q^k}\}_{k \in \Z}$ is an orthonormal basis for $\mathcal H$. Furthermore, for $z \in \{ z_-,z_+\}$ these functions satisfy the recurrence relation
\begin{gather} \label{eq:recrelvvbigqJac}
\big(x+x^{-1}\big) \psi_{zq^k}(x)=a_k \psi_{z q^{k+1}}(x) + b_k \psi_{z q^k}(x) + a_{k-1} \psi_{z q^{k-1}}(x),
\end{gather}
where
\begin{gather*}
a_k=a_k(a,c;z;q) = \left|\left(1-\frac{q^{-k}}{az}\right)\left(1-\frac{q^{-k}}{cz} \right) \right|,\\
b_k=b_k(a,c;z;q) = s^{-1}+s-s^{-1}\left|1-\frac{q^{1-k}}{az}\right|^2- s\left|1-\frac{q^{-k}}{cz}\right|^2.
\end{gather*}

We are now ready for the decomposition of $\mathcal T$. For this we need the vector-valued big $q^2$-Jacobi functions and the corresponding Hilbert space $\mathcal H(a,c;z_-,z_+|q^2)$, for certain values of the parameters $a$, $c$, $z_-$, $z_+$. Note that in this case the $q$ in all the formulas above must be replaced by $q^2$, e.g.~$s$ is given by $q|c|/|a|$.

As before, we f\/irst diagonalize $\mathcal T(\Omega)$ on a suitable subspace. For $p \in \Z$ we def\/ine the subspace
\[
\mathcal A_p = \mathrm{span}\{ f_{p,n}^- \mid n \in \Z\} \oplus \mathrm{span}\{ f_{p,n}^+ \mid n \in \Z\},
\]
where
\[
f_{p,n}^\sigma = (-1)^n e_{-n} \otimes e_{p+n}, \qquad n,p \in \Z,
\]
and we assume that $f_{p,n}^+$ is an element of the representation space of $\pi_{\rho_1,\eps_1+\eps}^{\mathrm P} \otimes \pi_{\rho_2,\eps_2-\eps}^{\mathrm P}$, and $f_{p,n}^-$ an element of the representation space of $\pi_{\rho_1',\eps_1-\eps}^{\mathrm P} \otimes \pi_{\rho_2',\eps_2+\eps}^{\mathrm P}$. Note that $\overline{\mathcal A_p} \cong \ell^2(\Z) \oplus \ell^2(\Z)$.

\begin{proposition} \label{prop:TePrincipal}
Let $a$, $c$, $z_-$, $z_+$ be given by
\begin{gather} \label{eq:acz-z+}
a=q^{2i\rho_2+2\eps_2+2p+1}, \qquad c=q^{2i\rho_1-2\eps_1+1}, \qquad z_\sigma= \sigma q^{-2\sigma \eps},  \qquad \sigma \in \{-,+\}.
\end{gather}
Then the operator $\Theta:\mathcal A_p \to \mathcal H(a,c;z_+,z_-|q^2)$ defined by
\[
f_{p,n}^\sigma \mapsto \big(x \mapsto \psi_{z_\sigma q^{2n}}\big(x;a,c|q^2\big)\big), \qquad \sigma \in \{-,+\},
\]
extends to a unitary operator. Furthermore, $\Theta$ intertwines $\mathcal T(\Omega)$ with $(q-q^{-1})^{-2} M_{x+x^{-1}-2}$.
\end{proposition}

Note in particular that the continuous spectrum occurs with multiplicity~2.
\begin{proof}
From \eqref{eq:De(OM)} and \eqref{eq:princ} we obtain
\[
\big(q-q^{-1}\big)^2 [\mathcal T(\Omega)+2 ]f_{p,n}^\sigma = a_n^\sigma f_{p,n+1}^\sigma + b_n^\sigma f_{p,n}^\sigma + a_{n-1}^\sigma f_{p,n-1}^\sigma,
\]
where
\begin{gather*}
a_n^\sigma = \left|(1-\sigma q^{-2n+2\eps_1+2\sigma\eps+2i\rho_1-1})(1-\sigma q^{-2n-2p-2\eps_2+2\sigma\eps+2i\rho_2-1})\right|,\\
b_n^\sigma  = q^{1-2\eps_1-2\eps_2-2p}+q^{-1+2\eps_1+2\eps_2+2p}-q^{1-2\eps_1-2\eps_2-2p}\left|1-\sigma q^{-2n+2\eps_1+2\sigma\eps+2i\rho_1-1} \right|^2\\
\phantom{b_n^\sigma  =}{} -q^{-1+2\eps_1+2\eps_2+2p}\left| 1-\sigma q^{-2n-2p+2\sigma\eps+2i\rho_2+1}\right|^2.
\end{gather*}
The result follows from comparing $\mathcal T(\Omega)f_{p,n}^\sigma$ with the recurrence relation \eqref{eq:recrelvvbigqJac}.
\end{proof}

The action of $K$ is given by
\[
\mathcal T(K) f_{p,n}^\sigma = q^{p+\eps_1+\eps_2} f_{p,n}^\sigma,
\]
and together with Proposition \ref{prop:TePrincipal} this leads to the following decomposition of $\mathcal T$.

\begin{theorem} \label{thm:decompTPP}
The $\U_q$-representation $\mathcal T_{\rho_1,\rho_2,\eps_1,\eps_2,\eps}$ is unitarily equivalent to
\[
2\dirint{0}{-\frac{\pi}{2 \ln q}} \pi^{\mathrm P}_{\rho,\eps'} d\rho\ \oplus\ \bigoplus_{\substack{j\in\Z\\ \sigma_j>0}} \pi_{\sigma_j,\eps'}^{\mathrm S} \ \oplus \bigoplus_{\substack{j\in\Z \\ k_j^+>\frac12}} \pi_{k_j^+}^+ \ \oplus \bigoplus_{\substack{j\in\Z \\ k_j^->\frac12}} \pi_{k_j^-}^-,
\]
where $\sigma_j=\eps_2-\eps_1+j-\frac12$, $\eps'=\eps_1+\eps_2$, $k_j^+=j+\eps_1+\eps_2$, $k_j^-=j-\eps_1-\eps_2$.
The unitary intertwiner is given by
\begin{gather*}
f_{p,n}^\sigma  \mapsto  \int_{0}^{-\frac{\pi}{2\ln q}}
\begin{pmatrix}
e_p^{\mathrm P} & 0\\
0 & e_p^{\mathrm P}
\end{pmatrix}
S\big(q^{2i\rho};a,c|q^2\big)\big(\Theta f_{p,n}^\sigma\big)\big(q^{2i\rho}\big) d\rho \\
\phantom{f_{p,n}^\sigma  \mapsto}{}
 + \sum_{\substack{j \in \Z\\ \sigma_j>0}} \big(\Theta f_{p,n}^\sigma\big)\big({-}q^{2\sigma_j}\big) \sqrt{w\big({-}q^{2\sigma_j};a,c;q^2\big)}   e_p^{\mathrm S} \\
\phantom{f_{p,n}^\sigma  \mapsto}{}
 + \sum_{\substack{j \in \N\\k_{p-j}^+>\frac12}} \big(\Theta f_{p,n}^\sigma \big)\big(q^{2k_{p-j}^+-1}\big) \sqrt{w\big(q^{2k_{p-j}^+-1};a,c;q^2\big)}  e_j^+ \\
\phantom{f_{p,n}^\sigma  \mapsto}{}
 + \sum_{\substack{j \in \N\\k_{-p-j}^->\frac12}} \big(\Theta f_{p,n}^\sigma \big)\big(q^{2k_{-p-j}^--1}\big) \sqrt{w\big(q^{2k_{-p-j}^--1};a,c;q^2\big)}  e_j^-
\end{gather*}
where $a$, $c$, $z_-$, $z_+$ are given by \eqref{eq:acz-z+}.
\end{theorem}

\begin{remark}
$(i)$ In this classical limit $q\uparrow 1$ the inf\/inite sum of strange series vanish, and we (formally) recover the decomposition of the tensor product of two principal unitary series of~$\su(1,1)$, see \cite[Theorem~II]{Pu} and \cite[Theorem~4.6]{Re}.

$(ii)$ Theorem \ref{thm:decompTPP} shows that the vector-valued big $q$-Jacobi functions have an interpretation as quantum analogs of Clebsch--Gordan coef\/f\/icients for~$\U_q$. In this interpretation the vector-valued big $q$-Jacobi function transform pair should be considered as a quantum analog of Neretin's~\cite{Ne} integral transform pair that has two $_3F_2$-functions as kernels.

$(iii)$ Note that the label $\eps$ does not appear in the decomposition of Theorem~\ref{thm:decompTPP}.
\end{remark}


\section{More quantum analogs of tensor product representations} \label{sec:more}
Using the (vector-valued) big $q$-Jacobi functions we can decompose several other quantum analogs of tensor product representations. We list a few decompositions here. The proofs are similar to the proofs of Theorems \ref{thm:decompT1} and \ref{thm:decompTPP}.

$\bullet$ A quantum analog of the tensor product of two complementary series:
for $\lambda_1,\lambda_2 \in (-\frac12,0)$, $\eps_1,\eps_2, \eps \in \R$,
\begin{gather*}
\Big(\big(\pi_{\lambda_1,\eps_1+\eps}^{\mathrm{C}}  \otimes \pi_{\lambda_2,\eps_2-\eps}^{\mathrm{C}}\big) \oplus \big(\pi_{\lambda_1+\frac12,\eps_1-\eps}^{\mathrm{S}} \otimes \pi_{\lambda_2+\frac12,\eps_2+\eps}^{\mathrm{S}}\big)\Big)\Delta
 \\
\qquad{} \cong 2\dirint{0}{-\frac{\pi}{2 \ln q}} \pi^{\mathrm P}_{\rho,\eps'} d\rho\ \oplus\ \bigoplus_{\substack{j\in\Z\\ \sigma_j>0}} \pi_{\sigma_j,\eps'}^{\mathrm S} \ \oplus \bigoplus_{\substack{j\in\Z \\ k_j^+>\frac12}} \pi_{k_j^+}^+ \ \oplus \bigoplus_{\substack{j\in\Z \\ k_j^->\frac12}} \pi_{k_j^-}^- \oplus \pi_{\lambda_1+\lambda_2,\eps'}^{\mathrm{C}},
\end{gather*}
under the condition $\lambda_1+\lambda_2<-\frac12$, where $\sigma_j=\eps_2-\eps_1+j-\frac12$, $\eps'=\eps_1+\eps_2$, $k_j^+=j+\eps_1+\eps_2$, $k_j^-=j-\eps_1-\eps_2$. If $\lambda_1+\lambda_2 \geq -\frac12$ the complementary series $\pi_{\lambda_1+\lambda_2,\eps'}^{\mathrm{C}}$ does not occur in the decomposition.

$\bullet$ A quantum analog of the tensor product of a principal unitary series and a complementary series:
for $\rho \in (0,-\frac{\pi}{2\ln q})$, $\lambda \in (-\frac12,0)$, $\eps_1,\eps_2,\eps \in \R$,
\begin{gather*}
\Big(\big(\pi_{\rho,\eps_1+\eps}^{\mathrm{P}}  \otimes \pi_{\lambda,\eps_2-\eps}^{\mathrm{C}}\big) \oplus \big(\pi_{\rho',\eps_1-\eps}^{\mathrm{P}} \otimes \pi_{\lambda+\frac12,\eps_2+\eps}^{\mathrm{S}}\big)\Big)\Delta
 \\
\qquad{} \cong 2\dirint{0}{-\frac{\pi}{2 \ln q}} \pi^{\mathrm P}_{\rho'',\eps'} d\rho''\ \oplus\ \bigoplus_{\substack{j\in\Z\\ \sigma_j>0}} \pi_{\sigma_j,\eps'}^{\mathrm S} \ \oplus \bigoplus_{\substack{j\in\Z \\ k_j^+>\frac12}} \pi_{k_j^+}^+ \ \oplus \bigoplus_{\substack{j\in\Z \\ k_j^->\frac12}} \pi_{k_j^-}^-,
\end{gather*}
where $\rho'=\rho-\frac{\pi}{2\ln q}$, $\sigma_j=\eps_2-\eps_1+j-\frac12$, $\eps'=\eps_1+\eps_2$, $k_j^+=j+\eps_1+\eps_2$, $k_j^-=j-\eps_1-\eps_2$.

$\bullet$ A quantum analog of the tensor product of a principal unitary series and a positive discrete series:
for $\rho \in (0,-\frac{\pi}{2\ln q})$, $k>0$, $\eps_1,\eps \in \R$,
\begin{gather*}
\Big(\big(\pi_{\rho,\eps_1}^{\mathrm{P}} \otimes \pi^+_k\big) \oplus \big(\pi_{\rho',\eps_1-\eps}^{\mathrm{P}} \otimes \pi_{k-\frac12,k+\eps}^{\mathrm{S}}\big)\Big)\Delta \cong
\dirint{0}{-\frac{\pi}{2 \ln q}} \pi^{\mathrm P}_{\rho'',\eps'} d\rho''\ \oplus\ \bigoplus_{\substack{j\in\Z\\ \sigma_j>0}} \pi_{\sigma_j,\eps'}^{\mathrm S} \oplus \bigoplus_{\substack{j\in\Z\\ k_j>\frac12}} \pi^+_{k_j}
\end{gather*}
where $\rho'=\rho-\frac{\pi}{2\ln q}$, $\sigma_j = k-\eps_1+\eps+j+\frac12$, $\eps'= \eps_1+k$, $k_j=k+\eps_1+j$.

$\bullet$ A quantum analog of the tensor product of a complementary series and a positive discrete series:
for $\lambda \in (-\frac12,0)$, $k>0$, $\eps_1,\eps \in \mathbb R$,
\begin{gather*}
\Big(\big(\pi_{\lambda,\eps_1}^{\mathrm{C}}  \otimes \pi^+_k\big) \oplus \big(\pi_{\lambda+\frac12,\eps_1-\eps}^{\mathrm{S}} \otimes \pi_{k-\frac12,k+\eps}^{\mathrm{S}}\big)\Big)\Delta \\
\qquad{} \cong  \dirint{0}{-\frac{\pi}{2 \ln q}} \pi^{\mathrm P}_{\rho,\eps'} d\rho\ \oplus\ \bigoplus_{\substack{j\in\Z\\ \sigma_j>0}} \pi_{\sigma_j,\eps'}^{\mathrm S} \oplus \bigoplus_{\substack{j\in\Z\\ k_j>\frac12}} \pi^+_{k_j} \oplus \pi_{-\lambda-k,\eps'}^{\mathrm{C}}
\end{gather*}
under the condition $\lambda+k<-\frac12$, where $\sigma_j = k-\eps_1+\eps+j+\frac12$, $\eps'= \eps_1+k$, $k_j=k+\eps_1+j$. If $\lambda + k \geq -\frac12$ the complementary series $\pi_{-\lambda-k,\eps'}^{\mathrm{C}}$ does not occur in the decomposition.

Let us remark that for the f\/irst two cases we need the vector-valued big $q$-Jacobi functions and corresponding orthogonality relations as described in Section~\ref{sec:PP}, but $a \neq \overline{b}$ and $c \neq \overline{d}$. In this case one extra discrete mass point can appear, and this occurs if $qd/as<1$, see~\cite{Gr09}. In the third case we need the big $q$-Jacobi functions and corresponding orthogonality relations as described in Section~\ref{sec:-+}, but with parameters $a$, $b$, $c$ satisfying $c=\overline{b}$, $a>0$ and $ab,ac,bc>1$. In this case the only discrete mass points are of the form $aq^k<1$, with $k \in \N$. Orthogonality relations for the big $q$-Jacobi functions with these conditions on $a$, $b$, $c$ can be obtained by very minor adjustments of the proof given in~\cite{KS03}.

\pdfbookmark[1]{References}{ref}
\LastPageEnding

\end{document}